\documentclass[a4paper,12pt]{article}
\usepackage[english]{babel}
\usepackage{amsmath,amsfonts,amssymb}

\usepackage[top=15mm,right=15mm,bottom=25mm,left=20mm]{geometry}

\newtheorem{theorem}{\hspace{\parindent}Theorem}
\newtheorem{lemma}{\hspace{\parindent}Lemma}
\newtheorem{corollary}{\hspace{\parindent}Corollary}[lemma]

\newtheorem{definition}{\hspace{\parindent}Definition}
\newenvironment{proof}[1][Proof]{
 \noindent\hspace{\parindent}\mbox{\textbf{#1.~}}}  {\hfill$\blacktriangleleft$}

\newcommand{\norm}[1]{\lVert#1\rVert}
\newcommand{\nnorm}[1]{\Bigl\|#1\Bigr\|}
\newcommand{\abs}[1]{\lvert#1\rvert}%
\newcommand{\angles}[1]{\langle{#1}\rangle}

\def\ot{\leftarrow}

\def\rr{\mathbf{R}}
\def\qq{\mathbf{Q}}
\def\zz{\mathbf{Z}}

\def\F{\mathcal{F}}
\def\sS{\mathcal{S}}
\def\E{\mathrm{E}}
\def\T{\mathbf{T}}
\def\F{\mathcal{F}}
\def\Z{\mathbf{Z}}
\def\R{\mathbf{R}}

\def\ffi{\varphi}

\begin{document}

\title{
 Fourier transform of function on locally compact Abelian groups
 taking value in Banach spaces
}
\author{
  Yauhen~Radyna, Anna~Sidorik\\
}
\date{}
\maketitle

\begin{abstract}
We consider Fourier transform of vector-valued functions on a
locally compact group~$G$, which take value in a Banach space $X$,
and are square-integrable in Bochner sense. If $G$ is a~finite
group then Fourier transform is a bounded operator. If $G$ is an
infinite group then  Fourier transform $F: L_2(G,X)\to
L_2(\widehat G,X)$ is a bounded operator if and only if Banach
space $X$ is isomorphic to a Hilbert one.
\end{abstract}

\section{Fourier transform over groups $\R$, $\Z$, $\T$}

In the paper \cite{Peetre} J.~Peetre proved an extension of
Hausdorf--Young's theorem describing image  of  $L_q(\rr)$ under
Fourier transform. He considered vector-valued  $x\in L_q(\rr,X)$,
$1\le q\le 2$ on the real axis taking value in Banach space $X$,
and integrable in Bochner sense, i.e. weakly measurable with
finite norm
\begin{equation*}
 ||x||_{L_q(\R,X)}=\Bigl(\int_{\rr}{||x(t)||_X^q d{t}}\Bigr)^{1/q}.
\end{equation*}
J.~Peetre noted, that for $q=2$ in all known to him cases Fourier
transform
\begin{equation*}
 \F:L_2(\rr,X)\to L_2(\rr,X), \quad
 (\F x)(s)= \int_{\rr}{x(t) e^{-2\pi i s t} d{t}}.
\end{equation*}
was bounded only if $X$ was isomorphic to a Hilbert space.
In~\cite{Kwapien} Polish mathematician S.~Kwapien in fact proved
the following

\begin{theorem}\label{th:Fourier-R}
Statements below are equivalent:

$1)$ Banach space $X$ is isomorphic to a Hilbert one.

$2)$ There exists $C>0$ such that for any positive integer $n$ and
$x_0, x_1,x_{-1}, \ldots,x_{n}, x_{-n} \in X$
\begin{equation*}
\int_0^1{\nnorm{\sum_{k=-n}^{n} e^{2\pi i k t}\cdot x_k}^2
d{t}}\le C\sum_{k = -n}^{n}\norm{x_k}^2.
\end{equation*}

$3)$ There exists $C>0$ such that for any positive integer $n$ and
$x_0, x_1,x_{-1}, \ldots,x_{n}, x_{-n} \in X$
\begin{equation*}
 C^{-1}\sum_{k = -n}^{n}\norm{x_k}^2
 \le
 \int_0^1{\nnorm{\sum_{k=-n}^{n} e^{2\pi i k t}\cdot x_k}^2 d{t}}.
\end{equation*}

$4)$ Fourier transform $\F$ initially defined on a dense subspace
of simple functions $D_{\F}\subset L_2(\rr,X)$,
\begin{equation*}
 D_{\F}= \Bigl\{ x(t)= \sum_{k=1}^{n}{I_{A_k}(t)\cdot x_k} \Bigr\},
\end{equation*}
is a bounded operator. Here $A_k$ are disjoint subset of finite
measure in $\rr$, $I_{A_k}$ are indicators (i.e. functions equal
to $1$ on $A_k$ and to $0$ elsewhere), $x_k$ are vectors in $X$.
\end{theorem}

Let us define Fourier transform of a vector-valued function over
integers $\zz$ by
\begin{equation*}
 \F_\Z: L_2(\Z,X)\to L_2(\T,X): (x_k) \mapsto \sum_{k\in\Z}{e^{-2\pi i k t}\cdot
 x_k},
\end{equation*}
Here $\T$ denotes a one-dimensional torus $\T=\R/\Z$, which is
isomorphic to $[0,1]$ with the length in the cathegory of spaces
with a measure.

Statement 2) in Theorem~\ref{th:Fourier-R} means that $\F_\Z
\mathcal{I}_\Z$ is a bounded operator $($here
$\mathcal{I}_\Z$~denotes the operator of changing variable
$\mathcal{I}_\Z: (x_k)\mapsto (x_{-k})$, which is an isometry$)$.
That is why $\F_Z$  is bounded on a dense subspace of $L_2(\Z,X)$
consisting of compactly-supported functions, and can be continued
to the whole $L_2(\Z,X)$.

Statement 3) in Theorem~\ref{th:Fourier-R} means that inverse
Fourier transform
\begin{equation*}
\F_\Z^{-1} : L_2(\T,X) \to L_2(\Z,X)
\end{equation*}
is a bounded operator.

\section{Generalization and necessary facts}

Fourier transform of Banach space-valued functions on a group
different from $\rr$,$\Z$,$\T$ (namely on additive group of
$p$-adic field $\qq_p$) was considered in~\cite{RRS-Fourier-Qp}.

Now it's natural to look at the general case of arbitrary locally
compact group $G$. We consider functions on $G$ taking value in
Banach space $X$, that are square-integrable in Bochner sense, and
Fourier transform
\begin{equation*}
 \F \equiv \F_G: L_2(G,X) \to L_2(\widehat{G},X),\quad
 (\F x)(\xi)=\int_{G}{\angles{\xi,t}_G x(t) d{\mu_G(t)}}.
\end{equation*}
Here $\widehat{G}$~is Pontryagin dual to $G$ (group of
characters), $\angles{\xi,t}_G$~is the canonical pairing between
$\widehat{G}$ and $G$, $\mu_G$~is Haar measure.

First we recall some necessary results. We refer to~\cite{H&R} for
results in harmonic analysis, and
to~\cite{Morris-DvoistvPontrygina} for structure theory of locally
compact groups, Bruhat--Schwartz theory is exposed
in~\cite{Bruhat-Distr}.

Fix a dual Haar measure $\mu_{\widehat{G}}$ on $\widehat{G}$ such
that scalar Steklov--Parseval's equality holds
\begin{equation*}
 ||\varphi||_{L_2(G)}^2= \int\limits_G{|\varphi|^2 d\mu_G}=
 \int\limits_{\widehat{G}}{|\F{\varphi}|^2 d\mu_{\widehat{G}}}= ||\F{\varphi}||_{L_2(\widehat{G})}.
\end{equation*}

Let $\sS(G)$ denote Bruhat--Schwartz space of ``smooth fastly
decreasing'' functions on $G$.

It's useful to take a dense subspace
\begin{equation*}
 D_\F= L_2(G)\otimes X \subset L_2(G,X),
\end{equation*}
as initial domain of $\F$, where it acts by
\begin{equation*}
 \F\Bigl(\sum_{k=1}^{n}{\ffi_k(t)\cdot x_k} \Bigr)=
 \sum_{k=1}^{n}{((\F\ffi_k)(\xi)\cdot x_k)}.
\end{equation*}
To show denseness of $D_\F$ and denseness of $\sS(G)\otimes X
\subset L_2(G)\otimes X$ consider indicator $I_A$ of arbitrary
measurable subset of finite measure (i.e. functions equal to $1$
on $A$ and to $0$ elsewhere). Clearly, $I_A\in L_2(G)$, and it can
be approximated by elements of $\sS(G)$ using convolution with any
$\delta$-net consisting of Bruhat--Schwartz functions. By
definition of $L_2(G,X)$ finite linear combinations
\begin{equation*}
 \sum_{k=1}^{n}{I_{A_k}(t)\cdot x_k}\in L_2(G)\otimes X  \equiv D_\F
\end{equation*}
are dense in $L_2(G,X)$. Here $A_k$~are measurable disjoint
subsets of finite measure in $G$, $I_{A_k}$ are indicators, and
$x_k\in X$.

Due to the fact that scalar Fourier transform is a bijection from
$L_2(G)$ into $L_2(\widehat{G})$, and is a bijection from $\sS(G)$
into $\sS(\widehat{G})$, the restriction of vector-valued Fourier
transform onto $L_2(G)\otimes X$ is a bijection into
$L_2(\widehat{G})\otimes X$, and its restriction onto
$\sS(G)\otimes X$ is a bijection to $\sS(\widehat{G})\otimes X$.

To handle the case of infinite group $G$ we need a theorem
describing  structure of locally compact Abelian
groups~\cite{Morris-DvoistvPontrygina}.

\begin{theorem}\label{th:osn-strukt}
Let $G$ be a locally compact Abelian group. Then  $G$ is a union
of open compactly generated subgroups $H$. Topology of $G$ is the
topology of inductive limit. In turn, each compactly-generated
subgroup $H\subset G$ is a projective limit of elementary
factor-groups $H/K$, where $K\subset H$ are compact.
``Elementary'' here means that $H/K$ is isomorphic to cartesian
product
\begin{equation*}
H/K \cong \mathbf{R}^{a_{H,K}} \times \mathbf{T}^{b_{H,K}} \times
\mathbf{Z}^{c_{H,K}}\times \mathbf{F}_{H,K},
\end{equation*}
where $a_{H,K}\ge 0$, $b_{H,K}\ge 0$, $c_{H,K}\ge 0$, and
$\mathbf{F}_{H,K}$ is a finite group.
\end{theorem}

\begin{definition}
We say that $G$ contains an $\mathbf{R}$-component if for some
$H,K$ in Theorem~\ref{th:osn-strukt} number $a_{H,K}$ is positive
in elementary factor-group  $\mathbf{R}^{a_{H,K}} \times
\mathbf{T}^{b_{H,K}} \times \mathbf{Z}^{c_{H,K}} \times
\mathbf{F}_{H,K}$.

In the same way we use phrases ``group $G$ contains
$\mathbf{Z}$-component'', ``group $G$ contains
$\mathbf{T}$-component''.
\end{definition}

Now recall some properties of Pontryagin duality.

Consider an exact sequence of homomorphisms (i.e. image of each
homomorphism coincides the kernel of the following one)
\begin{equation*}
  0\to K \to G\to G/K\to 0,\quad 0\to H \to G\to G/H\to 0,
\end{equation*}
where $K$ is compact, and $H$ is open. Dual groups
$\widehat{G/K}$, $\widehat{G/H}$ can be identified with
annihilators
\begin{equation*}
  K_G^\bot= \{ \chi\in\widehat{G}:\; \forall {g\in K},\; \angles{\chi,g}=1 \},\quad
  H_G^\bot= \{ \chi\in\widehat{G}:\; \forall {g\in H},\; \angles{\chi,g}=1 \}.
\end{equation*}
Moreover, $K_G^\bot$ is an open subgroup, and $H_G^\bot$ is
compact. One has dual oppositely-directed exact sequences
\begin{equation*}
  0\ot \widehat{G}/K_G^\bot \ot \widehat{G} \ot K_G^\bot \ot 0,\quad
  0\ot \widehat{G}/H_G^\bot \ot \widehat{G} \ot H_G^\bot \ot 0.
\end{equation*}
If $K\subset H$, then $K_G^\bot \supset H_G^\bot$.

It's not hard to see that Fourier transform of a function on $G$,
that is supported in open subgroup $H$ and is constant on cosets
of compact subgroup $K\subset H$, is a function on $\widehat{G}$
supported in $K_G^\bot$ and constant on cosets of $H_G^\bot$.

\section{The case of an arbitrary local compact group}

We will prove some necessary lemmas before formulating the main
result.

\begin{lemma}\label{lm:vec-Steklov-Parseval}
Assume that Banach space $X$ is isomorphic to a Hilbert one, i.e.
there exists inner product $(\cdot,\cdot)_X$ on $X$ such that for
some $C>0$ the following inequality is true
\begin{equation*}
 C^{-1}(x,x)_X^{1/2} \le ||x||_X \le C (x,x)_X^{1/2}.
\end{equation*}
Then Fourier transform $\F:L_2(G,X) \to L_2(\widehat{G},X)$ is
bounded.
\end{lemma}

\begin{proof}
Consider vector-valued Steklov--Parseval's equality
\begin{equation*}
  (\F \varphi, \F \varphi)_{L_2(\widehat{G},X)}=
  \int_{\widehat{G}}{(\F \varphi(\xi), \F \varphi(\xi)) d\mu_{\widehat{G}}(\xi)}=
  \int_{G}{(\varphi(t), \varphi(t)) d\mu_{G}(t)}=
  (\varphi,\varphi)_{L_2(G,X)},
\end{equation*}
which can be easily checked for $\varphi\in L_2(G)\otimes X$ by
means of axioms for inner product, scalar Steklov--Parseval's
equality and cross-norm's property
\begin{equation*}
  (\varphi_1\otimes x_1, \varphi_2\otimes x_2)_{L_2(G,X)}=
  (\varphi_1, \varphi_2)_{L_2(G)} \cdot   (x_1, x_2)_X.
\end{equation*}
Then we have on a dense subspace
\begin{equation*}
  ||\F\varphi||_{L_2(\widehat{G},X)} \le C (\F \varphi, \F \varphi)_{L_2(\widehat{G},X)}
  = C (\varphi, \varphi)_{L_2(G,X)} \le C^2
  ||\varphi||_{L_2(G,X)},
\end{equation*}
so we can extend $\F$ by continuity onto the whole $L_2(G,X)$ .
\end{proof}

\smallskip

If $G$ is a finite group, then space $L_2(G,X)$ is isomorphic to
the finite Cartesian product $X^G$. Pontryagin's dual group
$\widehat G$ is isomorphic to~$G$ itself. Self-dual Haar measure
possesses the property $\mu_G(G)=\sqrt{|G|}$. Fourier transform,
also known as discrete Fourier transform, becomes
\begin{equation*}
  (\F\varphi)(\xi)=
  \frac{1}{\sqrt{|G|}}{\sum_{t\in G} \angles{\xi,t}\varphi(t)}.
\end{equation*}

\begin{theorem}\label{th:F-finiteG}
If $G$~is a finite group, then Fourier transform $\F: L_2(G,X) \to
L_2(\widehat G,X)$ is bounded for any Banach space $X$.
\end{theorem}

\begin{proof}
It follows from inequality
\begin{gather*}
  ||\F \varphi||_{L_2(\widehat{G},X)}^2 =
  \sum_{\xi\in\widehat{G}}{\Bigl\|\frac{1}{\sqrt{|G|}}
     {\sum_{t\in G} \angles{\xi,t}\varphi(t)}\Bigr\|_X^2}\le
  \frac{1}{|G|}\sum_{\xi\in\widehat{G}}{\Bigl(
     \sum_{t\in G}{|\angles{\xi,t}|\cdot \|\varphi(t)\|_X}
  \Bigr) ^2}=
 \\
  = \frac{1}{|G|}\sum_{\xi\in\widehat{G}}{\Bigl(\sum_{t\in G}{\|\varphi(t)\|_X} \Bigr) ^2}
  = \Bigl(\sum_{t\in G}{\|\varphi(t)\|_X} \Bigr) ^2
  \le |G| \sum_{t\in G}{\|\varphi(t)\|_X^2} = |G|\cdot \| \varphi \|_{L_2(G,X)}^2.
\end{gather*}
\end{proof}

Now pass to infinite groups.

\begin{lemma}\label{lm:R-T-Z-sostav}
Let group $G$ contain  $\R$-component, $\T$-component or
$\Z$-component. Then boundedness of Fourier transform
\begin{equation*}
\F : L_2(G,X) \to L_2(\widehat{G},X)
\end{equation*}
implies isomorphism of Banach space $X$ to a Hilbert one.
\end{lemma}

\begin{proof}
Consider the case, when group $G$ contains  $\R$-component. There
are open compactly generated subgroup $H$ in $G$ and compact
subgroup $K\subset H$ such that $H/K\cong \R^a\times \T^b \times
\Z^c \times F$, where $a\ge 1$.

Let $\tau_1: H\to H/K$ be a canonical projection, $\tau_2: H/K\to
\R$ be the projection on the first coordinate of $\R^a$, $\tau=
\tau_2\circ\tau_1$.

Consider helper functions $\psi_{\R,i}\in L_2(\R)$, $2\le i\le
a$,\; $\psi_{\T,j}\in L_2(T)$, $1\le j\le b$,\; $\psi_{\Z,k}\in
L_2(\Z)$, $1\le k\le c$,\; $\psi_{F}\in L_2(F)$, each of them
having norm equal to $1$ in corresponding space. Consider
injection $J : L_2(\R,X)\to L_2(H,X)$,
\begin{equation*}
  J: \varphi \mapsto
  \left((\varphi\circ \tau) \otimes \Bigl(\bigotimes_{i=2}^{a}{\psi_{\R,i}}\Bigr)
  \otimes \Bigl( \bigotimes_{j=1}^{b}{\psi_{\T,j}} \Bigr)
  \otimes \Bigl( \bigotimes_{k=1}^{c}{\psi_{\Z,k}} \Bigr)
  \otimes \psi_{F} \right).
\end{equation*}
It is easy to see that the injection $J$ is isometric. There
exists a unique injection $\widehat{J}: L_2(\R,X)\to
L_2(\widehat{H},X)$, which is also an isometry, and for which the
following diagram is commutative
\begin{equation*}
  \begin{array}{ccc}
   L_2(\R,X)     & \overset{\mbox{$\F_{\R}$}}{\longrightarrow} &  L_2(\R,X) \\
   {J}\downarrow  &                         &  \downarrow {\widehat{J}} \\
   L_2(H,X)       & \overset{\mbox{$\F_{H}$}}{\longrightarrow} & L_2(\widehat{H},X).
  \end{array}
\end{equation*}

Space $L_2(H,X)$ can be identified with a closed subspace of
$L_2(G,X)$ consisting of functions, that are $0$ almost everywhere
outside $H$. Space $L_2(\widehat{H},X)$ can be identified with a
closed subspace of $L_2(\widehat{G},X)$ consisting of functions,
which are constant on cosets of $H^\bot$ (recall, that
$\widehat{H}\cong \widehat{G}/H_G^\bot$).

Fourier transform $\F_H$ is the restriction of $\F_G$ and thus,
bounded. Fourier transform $\F_{\R}=(\widehat{J})^{-1}\F_H J$ is
continuous as a composition of continuous maps. By
Theorem~\ref{th:Fourier-R} statement 4) space $X$ is isomorphic to
a Hilbert one.

Cases when group $G$ contains $\T$-component or $\Z$-component are
considered similarly.
\end{proof}

\bigskip

Consider the case, when group $G$ does not contain $\rr$-, $\zz$-
or $\T$-elements. Then all compactly generated subgroups $H\subset
G$ are projective limits of finite subgroups with discrete
topology, i.e. they are \emph{profinite} groups.

Profinite groups are characterized by the following
lemma~\cite{PRG}.

\begin{lemma}\label{lm:profinite}
Topological group $H$ is a profinite one if and only if it

a) possesses Hausdorff's property;

b) is compact;

c) is totally disconnected, i.e. for any two points $ x, y \in H$
there exists subset $U\subset H$ that is both open and closed,
such that $x \in U$ and $y \notin U$.
\end{lemma}

Any profinite group $H$ is either discrete (and finite by virtue
of compactness) or non-discrete (and therefore infinite).

Consider non-discrete profinite group $H$. We normalize Haar
measure on $H$ with $\mu_H(H)=1$.

Group $H$ is a Lebesgue space, i.e. it is isomorphic as a space
with measure to segment [0,1] with
length~\cite{H&R,Leonov-MatSociologia}. This fact can be proved
selecting sequence of compact subgroups $K_n\subset H$ such that
$K_1\subset K_2\subset ...$ and cardinality of quotient groups
$M_n:=|H/K_n|$ tends to $+\infty$. If one numbers cosets of $K_n$
properly, he becomes able to identify them  with the intervals of
length $1/M_n$ in $[0,1]$ up to a subsets of zero measure. By
$\tau$ denote this isomorphism of spaces with measure.

A system of functions $(r_i)_{i=1,2,...}$, similar to the
Rademacher's system on $[0,1]$ can be constructed on group $H$.
This is an orthogonal system of functions taking values
$\{+1,-1\}$ on subsets of measure~$\frac{1}{2}$. Saying in
probability-theoretical language functions $r_i$ are realizations
of independent random variables taking values $\{+1,-1\}$ with
probability $\frac{1}{2}$. One can simply assume $r_i=
r_i^\infty\circ\tau$, where $r_i^\infty(t)= \sin{2^i\pi t}$ on
$[0,1]$ is the usual Rademacher's functions.

We need a criterion, which is proved in~\cite{Kwapien}.

\begin{theorem}\label{th:E-criterium}
The following statements are equivalent:

1) Banach space $X$ is isomorphic to a Hilbert one.

2) There exists constant $C>0$ such that for any finite collection
of vectors $x_1, x_2, ..., x_n\in X$ two-sided Khinchin's type
inequality holds
\begin{equation*}
 C^{-1}\sum_{i=1}^{n}{\|x_i\|^2}\le
 \E\Bigl\|\sum_{i=1}^{n}{r_i x_i}\Bigr\|^2
 =\int_H\nnorm{\sum_{i = 1}^{n}r_i(t)x_i}^2 dt
 \le
 C\sum_{i=1}^{n}{\|x_i\|^2},
\end{equation*}
where $r_i$ are independent random variables taking values
$\{+1,-1\}$ with probability~$\frac{1}{2}$, and symbol $\E$
denotes expectation.
\end{theorem}

As in~\cite{Kwapien} we formulate a Lemma, which shows importance
of the system $(r_i)$ on $H$ and allows us to consider arbitrary
basis in $L_2(H)$ instead of $(r_i)$. By  $dt$ denote Haar measure
on~$H$.

\begin{lemma}\label{lm:r_i->f_i}
Let $X$ be a Banach space, $(f_i)$ be an orthonormal complete
system in $L_2(H)$. Assume that for some $C>0$ and for any
collection $x_1,x_2, \ldots,x_n \in X$ there is inequality
\begin{equation*}
 \int_H\nnorm{\sum_{i = 1}^{n}f_i(t)x_i}^2 dt
 \le C\sum_{i = 1}^{n}\norm{x_i}^2\quad
 \left(\text{resp., } C^{-1}\sum_{i = 1}^{n}\norm{x_i}^2 \le
  \int_H\nnorm{\sum_{i = 1}^{n}f_i(t)x_i}^2 dt \right).
\end{equation*}
Then for the same constant $C>0$ and for any collection $x_1,x_2,
\ldots,x_n \in X$ one also has
\begin{equation*}
 \int_H\nnorm{\sum_{i = 1}^{n}r_i(t)x_i}^2 dt \le
 C\sum_{i = 1}^{n}\norm{x_i}^2\quad
 \left(\text{resp., } C^{-1}\sum_{i = 1}^{n}\norm{x_i}^2 \le
  \int_H\nnorm{\sum_{i = 1}^{n}f_i(t)x_i}^2 dt \right).
\end{equation*}
\end{lemma}

\begin{proof}
As Rademacher's system $(r_i)$ is orthonormal and $(f_k)$ are
complete, we can find for a given $\varepsilon>0$ an increasing
sequences of indices $(k_j)$, $(m_j)$ and orthonormal sequence
$(h_j)$ such that
\begin{equation*}
h_j=\sum_{k=k_j}^{k_{j+1}-1}(h_j,f_k)\cdot f_k,\    \
\int_{H}\abs{h_j(t)-r_{m_j}(t)}^2dt<\frac{\varepsilon}{2^{j}}.
\end{equation*}
For a fixed $n$ and fixed $x_1,x_2, \ldots,x_n \in X$ we have
\begin{equation*}
\int_{H}{\nnorm{\sum_{i=1}^{n}r_i(t)\cdot x_i}^2 d{t}} =
\int_{H}{\nnorm{\sum_{i=1}^{n}r_{m_j}(t)\cdot x_i}^2 d{t}}.
\end{equation*}
By the triangle inequality
\begin{gather*}
 \Bigl(\int_{H}\nnorm{\sum_{j=1}^{n}r_{m_j}(t)\cdot x_j}^2 d{t}\Bigr)^{1/2}\le
 \Bigl(\int_{H}\nnorm{\sum_{j=1}^{n}(r_{m_j}(t)-h_j(t))\cdot x_j}^2
 d{t}\Bigr)^{1/2} +
\\
 + \Bigl(\int_{H}\nnorm{\sum_{j=1}^{n}h_{j}(t)\cdot x_j}^2 d{t}\Bigr)^{1/2}
 \le \sqrt{\varepsilon}\Bigl(\sum_{j=1}^{n}\norm{x_j}^2\Bigl)^{1/2}
 +\Bigl(\int_{H}\nnorm{\sum_{j=1}^{n}h_{j}(t)\cdot x_j}^2 d{t}\Bigr)^{1/2}.
\end{gather*}
As $\displaystyle
1=\norm{h_j}^2=\sum_{k=k_j}^{k_{j+1}-1}\abs{(h_j,f_k)}^2$, we get
\begin{gather*}
 \int_{H}\nnorm{\sum_{j = 1}^{n}h_{j}(t)\cdot x_j}^2dt=
 \int_{H}\nnorm{\sum_{j = 1}^{n}
 \biggl(\sum_{k=k_j}^{k_{j+1}-1}(h_j,f_k)\cdot f_k\biggr)x_j}^2 d{t}
 \le \\ \le
 C \sum_{j=1}^{n}\sum_{k=k_j}^{k_{j+1}-1}\abs{(h_j,f_k)}^2\norm{x_j}^2=
 C \sum_{j=1}^{n}\norm{x_j}^2.
\end{gather*}
Thus,
\begin{equation*}
 \int_{H}{\nnorm{\sum_{i = 1}^{n}r_i(t)\cdot x_i}^2 d{t}}\le
 (\sqrt{\varepsilon}+\sqrt{C})^2\sum_{i = 1}^{n}\norm{x_i}^2.
\end{equation*}
As $\varepsilon$ is arbitrary, we obtain the desired inequality.
Proof in the case of reverse type inequality is analogous.
\end{proof}

\begin{corollary}\label{co:ner}
Let $X$ be a Banach space, and $(f_i)$ be a complete orthonormal
system in $L_2(H)$. Space $X$ is linearly isomorphic to a Hilbert
one if and only if there exists constant $C>0$ such that for any
set of vectors $x_1,x_2, \ldots,x_n \in X$ one has
\begin{equation*}\label{4}
 C^{-1}\sum_{i = 1}^{n}\norm{x_i}^2\le
 \int_{H}\nnorm{\sum_{i=1}^{n}f_i(t)\cdot x_i}^2 dt
 \le C\sum_{i = 1}^{n}\norm{x_i}^2.
\end{equation*}
\end{corollary}

In Corollary~\ref{co:ner} isomorphism of $X$ to a Hilbert space
follows from \emph{two-sided} inequality. Knowledge of profinite
groups' structure allows us to switch from lower estumate to upper
one and vice versa as shown below.

\bigskip

First recall that Bruhat--Schwartz space on a profinite group $H$
and on dual discrete $\widehat{H}$ consists of locally constant
functions with compact support. Spaces $\sS(H)$ and
$\sS(\widehat{H})$ are inductive limit of finite-dimensional
spaces and carry the strongest locally convex
topology~\cite{Bruhat-Distr}.

\bigskip

Now we are going to study the case of vector-valued Fourier
transform on a profinite non-discrete group~$H$. As $H$ is a
compact infinite group, its dual $\widehat H$ is a discrete
infinite group.

\begin{lemma}\label{lm:L_2(H,X)-criterium}
Let $X$ be a Banach space, and $H$ be a profinite non-discrete
group. The following statements are equivalent:

$1)$ $X$ is linearly isomorphic to a Hilbert space.

$2)$ There exists constant $C>0$ such that for any set of vectors
$x_1,...,x_n\in X$ and characters $\xi_1,...,\xi_n\in \widehat{H}$
one has
\[
 \int_{H}\nnorm{\sum_{k = 1}^{n} \langle \xi_k,t \rangle x_k}^2 dt
 \le C \sum_{k = 1}^{n}\norm{x_k}^2.
\]

$2)'$ Fourier transform $\F_{\widehat{H}}: L_2(\widehat{H},X)\to
L_2(H,X)$ and inverse Fourier transform $F_{H}^{-1}=\mathcal{I}_H
\F_{\widehat{H}}$ are bounded. Here $\mathcal{I}_H$ is an
isometrical operator of changing variable $x(t)\mapsto x(-t)$.

$3)$ There exists constat $C>0$ such that for any set of vectors
$x_1,...,x_n\in X$ and characters  $\xi_1,...,\xi_n\in
\widehat{H}$ one has
\[
 C^{-1} \sum_{k = 1}^{n}\norm{x_k}^2 \le
 \int_{H}\nnorm{\sum_{k = 1}^{n} \langle \xi_k,t \rangle x_k}^2
 dt.
\]

$3)'$ Inverse Fourier transform $\F_{\widehat{H}}^{-1} : L_2(H,X)
\to L_2(\widehat{H},X)$ and Fourier transform $\F_H= \mathcal{I}_H
\F_{\widehat{H}}^{-1}$ are bounded.
\end{lemma}

\begin{proof}
Lemma~\ref{lm:vec-Steklov-Parseval} yields implications
$1)\Rightarrow 2)'$, $1)\Rightarrow 3)'$.

To prove equivalences  $2)'\Leftrightarrow 2)$, $3)'\Rightarrow
3)$ it's enough to consider
\begin{equation*}
h= \sum_{k=1}^{n}{ I_{\{\xi_k\}}\cdot x_k } \in
\sS(\widehat{H})\otimes X,
\end{equation*}
where $\xi_k \in \widehat{H},  x_k \in X, n \in \mathbf{N}$. By
definition of Fourier transform $\F: L_2(\widehat{H},X)\to
L_2(H,X)$ and by chosen normalization of Haar measure on $H$ we
get
\begin{equation*}
 \norm{h}^2_{L_2(\widehat{H},X)}=
 \sum_{k = 1}^{n}\norm{x_k}^2,\quad \norm{\F h}^2_{L_2(H,X)}=
 \int_{H}\nnorm{\sum_{k = 1}^{n} \langle \xi_k,t \rangle x_k}^2 dt.
\end{equation*}
Equivalence follows from the density of $\sS(\widehat{H})\otimes
X$ in $L_2(\widehat{H},X)$, density of $\sS(H)\otimes X$ in
$L_2(H,X)$, and  bijectivity of  $\F_{\widehat{H}}$, $\F_H$ in the
corresponding spaces.

By Corollary~\ref{co:ner} one also has $2)\;\&\;3) \Rightarrow
1)$.

Now assume boundedness condition $2)'$. Let's look at Fourier
transform $\F_H$ on a subspace $\sS(H)\otimes X$. For this
consider arbitrary compact subgroup $K\subset H$, for which
$|H/K|<+\infty$. We identify functions, which are constant on the
cosets of $K$, with elements of $\sS(H/K)\otimes X$.

By finiteness of $H/K$ there exist an isomorphism
$\alpha:\widehat{H/K}\to H/K$. Adjoint isomorphism
$\alpha^*:\widehat{H/K}\to H/K$ is defined by
\[
 \angles{\xi_1,\alpha(\xi_2)}_H=
 \angles{\alpha^*(\xi_1),\xi_2}_{\widehat{H}}.
\]

Consider operator
\[
 R_\alpha: \sS(H/K)\otimes X\to \sS(\widehat{H/K})\otimes X:
 (R_\alpha\psi)(\xi')= \psi(\alpha(\xi'))\cdot |H/K|^{-\frac{1}{2}}.
\]
It's an isometry, because
\begin{gather*}
  ||R_\alpha\psi||^2 =
  \sum_{\xi'\in\widehat{H/K}}{||R_\alpha\psi(\xi)||^2
  \mu_{\widehat{H/K}}(\xi)}=
  \sum_{\xi'\in\widehat{H/K}}{||\psi(\alpha(\xi'))||^2\cdot
  |H/K|^{-\frac{1}{2}\cdot 2}}=
\\=
 [t:=\alpha(\xi')] =
 \sum_{t\in{H/K}}{||\psi(t)||^2 \cdot |H/K|^{-1}}=
 \sum_{t\in{H/K}}{||\psi(t)||^2 \mu_{H/K}(t)}=
 ||\psi||^2.
\end{gather*}

Now one has
\begin{gather*}
 (\F_{\widehat{H}} R_\alpha \psi)(t')=
 \sum_{\xi'\in\widehat{H/K}}{\angles{t',\xi'}_{\widehat{H/K}}(\psi(\alpha(\xi'))|H/K|^{-\frac{1}{2}})},
 \\
 (R_{(\alpha^*)}\F_{\widehat{H}} R_\alpha \psi)(\xi)=
 \Bigl(
    \sum_{\xi'\in\widehat{H/K}}
    {\angles{\alpha^*(\xi),\xi'}_{\widehat{H/K}}\psi(\alpha(\xi'))|H/K|^{-\frac{1}{2}}}
    \Bigr) |H/K|^{-\frac{1}{2}}
 = [ t: = \alpha(\xi')] = \\=
 \Bigl(
    \sum_{t\in{H/K}}
    {\angles{\alpha^*(\xi),\alpha^{-1}(t)}_{\widehat{H/K}}\psi(t) } \Bigr) |H/K|^{-1}=
    \sum_{t\in{H/K}}
    {\angles{\alpha^*(\xi),\alpha(\alpha^{-1}(t))}_{H/K}\psi(t) } |H/K|^{-1}=
 \\ =
    \sum_{t\in{H/K}}
    {\angles{\xi,t}_{H/K}\psi(t)\mu_{H/K}(t) } =
    (\F_H \psi)(\xi).
\end{gather*}

Thus, restriction of $\F_H$ onto $\sS(H/K)\otimes X$ has the same
norm as $\F_{\widehat{H}}$ does. As $K$ is arbitrary, it implies
continuity of $\F_H$ on $L_2(H,X)$, and implication $2)\Rightarrow
3)$ is true. Implication $3)\Rightarrow 2)$ can be proved in the
same way.

Now one has enough implications to see the equivalence of all
statements.
\end{proof}

\bigskip

Let's pass to the general case of Fourier transform on a locally
compact Abelian group $G$. This is our main result.

\begin{theorem}
Let  $X$ be a Banach space, let $G$ be a locally compact Abelian
group. Space $X$ is linearly isomorphic to a Hilbert one if and
only if Fourier transform
\begin{equation*}
\F : L_2(G,X) \to L_2(\widehat G,X)
\end{equation*}
is bounded.
\end{theorem}

\begin{proof}
Existence of isomorphism is a sufficient condition for boundedness
of $\F$ as shown in Lemma~\ref{lm:vec-Steklov-Parseval}.

Now assume that Fourier transform is bounded.

If $G$ contains \mbox{$\rr$-,} $\T$- or $\zz$-component, then
isomorphism of $X$ to a Hilbert space follows from
Lemma~\ref{lm:R-T-Z-sostav}, and we are done.

Otherwise, group $G$ does not contain \mbox{$\rr\-$,} $\T\-$ or
$\zz$-components. In~this case all compactly generated open
subgroups $H\subset G$ are profinite.

If some of these $H$ is non-discrete, then we consider space
$L_2(H,X)$ as a closed subspace in $L_2(G,X)$ (one can simply
assume that function from $L_2(H,X)$ are zero outside $H$). We
identify $L_2(\widehat{H},X)$ with a subspace of
$L_2(\widehat{G},X)$ consisting of functions constant on cosets of
annihilator $H_G^\bot\subset\widehat{G}$. Fourier transform on
$L_2(H,X)$ is the restriction of Fourier transform from $L_2(G,X)$
and is bounded. Isomorphism of $X$ to a Hilbert space follows from
Lemma~\ref{lm:L_2(H,X)-criterium}, statement $3)'$. And we are
done.

If all subgroups $H\subset G$ considered are discrete (their
compactness implies finiteness), then by
Theorem~\ref{th:osn-strukt} group $G$ is an \emph{inductive} limit
of discrete subgroups. By properties of Pontryagin duality (in the
language of Cathegory theory one can say that passing to a dual
group is an exact functor) dual group  $\widehat{G}$ is a
\emph{projective} limit of $\widehat{H}$. Groups $\widehat{H}$ are
dual to finite discrete $H$. Thus, $\widehat{H}$ are isomorphic
$H$, and are finite discrete themselves. Group $\widehat{G}$ is
profinite. As $G$ is infinite, $\widehat{G}$ is non-discrete.

Boundedness of Fourier transform $\F: L_2(G,X)\to
L_2(\widehat{G},X)$ is equivalent to boundedness of inverse
Fourier transform  $\F_{\widehat{G}}^{-1}$ on profinite
non-discrete $\widehat{G}$. Isomorphism of $X$ to a Hilbert space
follows from Lemma~\ref{lm:L_2(H,X)-criterium} statement $2)'$.
\end{proof}

\end{document}